\documentclass[12pt,oneside,reqno]{amsart}
\usepackage{amssymb}
\newtheorem{thm}{Theorem}
\newtheorem{lem}[thm]{Lemma}
\newtheorem{pro}[thm]{Proposition}

\newtheorem{con}[thm]{Conjecture}
\theoremstyle{remark}
\newtheorem{remark}[thm]{Remark}

\begin{document}
\title[Oppenheim conjecture for pairs]{Oppenheim conjecture for pairs consisting of a linear form and a quadratic form}
\author{Alexander Gorodnik}
\address{Department of Mathematics, University of Michigan, Ann Arbor, MI 48109}
\email{gorodnik@umich.edu}
\thanks{This article is a part of author's PhD thesis at Ohio State University done under supervision of Prof.~Bergelson.}

\begin{abstract}
Let $Q$ be a nondegenerate quadratic form, and $L$ is a nonzero linear form of dimension $d>3$.
As a generalization of the Oppenheim conjecture, we prove that the set $\{(Q(x),L(x)):x\in\mathbb{Z}^d\}$
is dense in $\mathbb{R}^2$ provided that $Q$ and $L$ satisfy some natural conditions.
The proof uses dynamics on homogeneous spaces of Lie groups.
\end{abstract}

\maketitle

\section{Introduction}

It was proved by Margulis \cite{m89} (see \cite{mar97} for an up-to-date survey) that 
if $Q$ is a real indefinite nondegenerate quadratic form of dimension $d\ge 3$ which is
not proportional to a rational form, then $\{Q(x):x\in\mathbb{Z}^d\}$ is dense in $\mathbb{R}$.
A similar problem was considered for pairs $(Q,L)$ where $Q$ is a quadratic form as above,
and $L$ is a nonzero linear form. The known results are limited to dimension $3$.
It was proved by Dani and Margulis \cite{dm90} that $\{(Q(x),L(x)):x\in\mathbb{Z}^3\}$
is dense in $\mathbb{R}^2$ if no nonzero linear combination of $Q$ and $L^2$ is rational, and
the plane $\{L=0\}$ is tangent to the surface $\{Q=0\}$. Clearly, the first condition is neccesary for the density to hold.
The second condition guarantees that the stabilizer in $\hbox{SL}(3,\mathbb{R})$
of the pair $(Q,L)$ is unipotent, so that the results on orbits of unipotent flows can be used.
One can hope to remove the second condition. However, Dani proved in \cite{dan00} that
if the surface $\{Q=0\}$ and the plane $\{L=0\}$ intersect transversally, the density can fail for a set pairs of
full Hausdorff dimension.
On the other hand, it is easy to see using Moore ergodicity criterion
that the density holds for a set of pairs of full measure provided that the surfaces
$\{Q=0\}$ and $\{L=0\}$ have nonzero intersection.

Denote by $\mathcal{P}(\mathbb{Z}^d)$ the set of primitive integer vectors in $\mathbb{Z}^d$.
The results mentioned above still hold when $\mathbb{Z}^d$ is replaced by $\mathcal{P}(\mathbb{Z}^d)$.

In this paper we prove a density result for pairs consisting of a linear form and a quadratic form
of dimension $d \ge 4$:

\begin{thm} \label{ql_m}
Let $Q$ be a nondegenerate quadratic form, and $L$ be a nonzero linear form in
dimension $d\ge 4$. Suppose that
\begin{enumerate}
\item[1.] $Q|_{L=0}$ is indefinite.
\item[2.] 
For every $(\alpha,\beta)\ne (0,0)$, with $\alpha,\beta\in\mathbb{R}$,
$\alpha Q+\beta L^2$ is not rational.
\end{enumerate}
Then $\{(Q(x),L(x)):x\in\mathcal{P}(\mathbb{Z}^d\})$ is dense in $\mathbb{R}^2$.
\end{thm}

For $d=3$, condition (1) implies that the surfaces $\{Q=0\}$ and $\{L=0\}$ intersect transversally.
Therefore, it follows from the result of Dani \cite{dan00} mentioned above that 
the analogue of Theorem \ref{ql_m} does not hold for $d=3$.

It is easy to see that the second condition is neccesary for the conclusion of the theorem to hold.
A condition similar to the first condition is required to insure that 
$$
\{(Q(x),L(x)): x\in \mathbb{R}^d\}= \mathbb{R}^2.
$$
It is possible that the first condition can be weakened (see Conjecture \ref{ql_con} below).

{\sc Acknowledgment:} The author is deeply grateful to his advisor V.~Bergelson for suggesting the
problem as well as for many helpful discussions.

\section{Canonical forms}

Let $Q_i$, $i=1,2$, be quadratic forms, and  $L_i$, $i=1,2$, be linear forms of dimension $d$.
We say that the pairs $(Q_1,L_1)$ and $(Q_2,L_2)$ are {\it equivalent}
if $Q_1(x)=\lambda Q_2(g\cdot x)$ and $L_1(x)=\mu L_2(g\cdot x)$ for some
$\lambda, \mu\in\mathbb{R}-\{0\}$ and $g\in\hbox{SL}(n,\mathbb{R})$. That is, they can be transformed
into each other by a linear change of coordinates and scaling.

\begin{pro} \label{ql_can}
Every pair $(Q,L)$, where $Q$ is a nondegenerate quadratic form, and $L$ is a nonzero linear form, is equivalent
to one and only one of the following pairs:
\begin{enumerate}
\item[(I)] $\left(x_1^2+\ldots +x_s^2-x_{s+1}^2-\ldots-x_d^2, x_d\right)$\\ where $s=1,\ldots,d$.
\item[(II)] 
$\left(x_1^2+\ldots +x_s^2-x_{s+1}^2-\ldots-x_{d-2}^2+x_{d-1}x_d, x_d\right)$ where $s=0,\ldots,[\frac{d-2}{2}]$.
\end{enumerate}
\end{pro}

\begin{proof}
Without loss of generality, $L(x)=x_d$. Applying a linear transformation
in $x_1,\ldots,x_{d-1}$, we transform $Q$ to the form 
$$
\sum_{i\le k} \pm x_i^2 + \left(\sum_j a_j x_j\right)x_d\;\;\hbox{for some}\;\; k\le d-1.
$$
Using linear transformations $x_i\rightarrow x_i\pm\frac{a_i}{2}x_d$ for $i\le k$, we can make $a_i=0$ for $i\le k$.
If $a_i=0$ for all $i\ne d$, then $k=d-1$ because $Q$ is nondegenerate. Thus, we are in the case (I).
Let $a_l\ne 0$ for some $l<d$. Making a linear change of variables $x_l\leftarrow\sum_j a_j x_j$, we get
$$
Q=\sum_{i\le k} \pm x_i^2 + x_lx_d.
$$
Since $Q$ is nondegenerate, $k=d-2$ and $l=d-1$.
This is the case (II). Note that a pair $(Q,L)$ is of type (I) iff $Q+\alpha L^2$ is degenerate for some
$\alpha\in\mathbb{R}$. In particular, this shows that pairs of type (I) and (II) are not equivalent.
Signatures of quadratic forms $Q$ and $Q|_{L=0}$ are invariants of equivalence
with one possible exception that pairs with signatures $(s,d-s)$ and $(d-s,s)$ could be equivalent.
For example, pairs of type (I) with $s=0$ and $s=d$ are equivalent.
Let $(Q,L)$ be as in (I) with $s=i$ for $i\ne 0,d$. Then $Q|_{L=0}$ has signature $(i,d-i-1)$.
On the other hand, when $s=d-i$, $Q|_{L=0}$ has signature $(d-i,i-1)$.
Clearly, these two cases are not equivalent unless $n=2i$, and they coincide.
When $(Q,L)$ is of type (II), it is easy to see that the cases $s=i$ and $s=d-2-i$ are equivalent
for $i=0,\ldots, d-2$. This finishes the proof.
\end{proof}

%
%
%

Using the explicit formulas for canonical forms, it is easy to prove
the following lemma:

\begin{lem} \label{ql_im}
Let $Q$ be a nondegenerate quadratic form, and $L$ is a nonzero linear form
such that $Q|_{L=0}$ is indefinite. Then
$$
\{(Q(x),L(x)):x\in\mathbb{R}^d\}=\mathbb{R}^2.
$$
\end{lem}




\section{Reduction to lower dimension}

In this section, we prove Proposition \ref{ql_red}, which reduces the proof of Theorem \ref{ql_m}
to the case of dimension $4$.

\begin{pro} \label{ql_red}
Let $(Q,L)$ satisfy the conditions of Theorem \ref{ql_m}.
If the dimension of the space is greater than $4$, then there exists a rational subspace $U$
of codimension $1$ such that $\left(Q|_U, L|_U\right)$ satisfies the same conditions.
Moreover, $U$ can be chosen such that $\left(Q|_U, L|_U\right)$ is of type (I).
\end{pro}

Consider the set $\Omega=\{(\ell_1,\ell_2,\ell_3)\}\subset\mathbb{R}^{3d}$ of triples of linear forms
of dimension $d$ defined by the following conditions:
\begin{enumerate}
\item[1.] $\ell_1,\ell_2,\ell_3$ are linearly independent.
\item[2.] $L\ne 0$ on the space $\{\ell_i=0:i=1,2,3\}$.
\item[3.] $Q$ is nondegenerate on the space $\{\ell=0:\ell\in T\}$ for any subset $T\subseteq\{\ell_1,\ell_2,\ell_3\}$.
\item[4.] $\left(Q|_{\ell_i=0},L|_{\ell_i=0}\right)$ is of type (I) for $i=1,2,3$.
\item[5.] $Q|_{\{\ell_i=0\}\cap\{L=0\}}$ is indefinite for $i=1,2,3$.
\end{enumerate}

\begin{lem}
The set $\Omega$ is nonempty and open.
\end{lem}

\begin{proof}
Let $\Omega_0$ be the set of triples of linear forms that satisfy conditions (1) and (2), and
$\Omega_j$ the set of triples in $\Omega_0$ that satisfy condition ($j$), $j=3,4,5$.
Note that the set $\Omega_0$ is the complement of a proper algebraic set
(i.e. a set defined by polynomial equations).

The sets $\Omega_j$, $j=3,4,5$, are not empty. This is easy to see using the canonical forms from
Proposition \ref{ql_can}. For example,
when $(Q,L)$ is the canonical form (II), $(\ell_1,\ell_2,\ell_3)\in\Omega_4$ for
$\ell_i(\bar x)=x_{d-4-i}+x_{d-1}+x_d$, $i=1,2,3$.
It is also clear that the sets $\Omega_j$, $j=3,4,5$, are open.
To show this for $\Omega_4$, we use that 
a pair $(Q,L)$ is of type (I) iff
$$
\det (Q)\ne 0,\; \det(Q|_{L=0})\ne 0,\; L\ne 0.
$$
It suffices to show that $\Omega_3$ and $\Omega_4$ are dense.
Indeed, then it follows that $\Omega_3\cap\Omega_4$
is an open dense set, and $\Omega=\Omega_3\cap\Omega_4\cap\Omega_5$ is nonempty and open.

We claim that for every $T\subseteq\{\ell_1,\ell_2,\ell_3\}$, the set of triples $(\ell_1,\ell_2,\ell_3)\in\Omega_0$
such that $Q$ is nondegenerate on the space $\{\ell=0:\ell\in T\}$ contains the complement of
a proper algebraic subset. This will imply that the set $\Omega_3$ is dense.
To simplify notations, we only consider the case of one linear form.
Namely, we show that the set of nonzero linear forms $\ell$ such that $Q|_{\ell=0}$ is
nondegenerate contains the complement of a proper algebraic subset. The cases of pairs and triples are handled similarly.
For a nonzero linear form $\ell$, there exists a basis $\{v_i^{\ell}:i=1,\ldots,d-1\}$ of the space 
$\{\ell=0\}$ such that the coordinates of the vectors $v_i^{\ell}$ are rational
functions of the coefficients of $\ell$. Define a rational function
$$
\phi(\ell)\stackrel{def}{=}\det\left(Q(v_i^{\ell},v_j^{\ell}):i,j=1,\ldots,d-1\right).
$$
For a linear form $\ell$ such that $\phi(\ell)$ is defined, $Q|_{\ell=0}$ is nondegenerate iff
$\phi(\ell)\ne 0$.
If $\phi=0$ on its domain, $Q|_{\ell=0}$ is nondegenerate only for $\ell$ in a proper
algebraic subset where $\phi$ is undefined.
This is a contradiction because $\Omega_3$ is nonempty and open. Hence, $\phi\ne 0$.
The form $Q|_{\ell=0}$ is nondegenerate for $\ell\in\{\ell:\phi(\ell)\ne 0\}$. This set is
the complement of a proper algebraic set.

To prove that $\Omega_4$ is dense, we show that the set of nonzero linear forms $\ell$ such that
$(Q|_{\ell=0},L|_{\ell=0})$ is of type (I) contains the complement of a proper algebraic subset.
The form $(Q|_{\ell=0},L|_{\ell=0})$ is of type (I) iff
$$
\det (Q|_{\ell=0})\ne 0,\; \det(Q|_{\{L=0,\ell=0\}})\ne 0,\; L|_{\ell=0}\ne 0.
$$
Note that the first inequality and the third inequality define sets
that contain the complements of proper algebraic subsets. Hence,
it is enough to prove that the set of linear forms $\ell$ that satisfy the second inequality
contains the complement of a proper algebraic subset. 
There exists a basis $\{w_i^{\ell}:i=1,\ldots,d-2\}$ of the space $\{L=0,\ell=0\}$ such that
the coordinates of the basis vectors are rational functions of the coefficients of $\ell$.
Define a rational function 
$$
\psi(\ell)\stackrel{def}{=}\det\left(Q(w_i^{\ell},w_j^{\ell}):i,j=1,\ldots,d-2\right).
$$
For $\ell\;$ such that $\psi(\ell)$ is defined, $Q|_{\{L=0,\ell=0\}}$ is nondegenerate
iff $\psi(\ell)\ne 0$. Therefore, we deduce as above that the set of $\ell$ such that
$Q|_{\{L=0,\ell=0\}}$ is nondegenerate contains the complement of a proper algebraic
subset. This implies that $\Omega_4$ is dense and proves the lemma.
\end{proof}

Denote by $\pi:\mathbb{R}^d\times\mathbb{R}^d\times\mathbb{R}^d\rightarrow\mathbb{R}^d\times\mathbb{R}^d$
the projection on the first two coordinates, and
by $p:\mathbb{R}^d\times\mathbb{R}^d\times\mathbb{R}^d\rightarrow\mathbb{R}^d$
the projection on the first coordinate.

\begin{lem} \label{ql_red1}
Suppose that $(l_1,m_1),(l_2,m_2),(m_1,m_2)\in\pi (\Omega)$.
Then there exists a rational linear form $n$ such that 
$$
(l_1,m_1,n),(l_2,m_2,n),(m_1,m_2,n)\in\Omega.
$$
\end{lem}

\begin{proof}
As in the previous lemma, one shows that if $(l,m)\in\pi (\Omega)$, the set of $n$
such that $(l,m,n)$ satisfies conditions (1)--(4) is open and dense. Thus, the set of
linear forms $n$ such that the triples $(l_1,m_1,n)$, $(l_2,m_2,n)$, and $(m_1,m_2,n)$
satisfy conditions (1)--(4) is open and dense. The condition that $Q|_{\{n=0\}\cap\{L=0\}}$
is indefinite holds on an open set of $n$. Hence, such a linear form $n$ exists.
\end{proof}

For a linear form $\ell$, consider a linear map $F_\ell$ from the space
of quadratic forms on $\mathbb{R}^d$ to the space of quadratic forms on $\{\ell=0\}$
defined by 
$$
F_\ell(Q)=Q|_{\ell=0}.
$$
This map is rational if $\ell$ is
rational. The kernel of $F_\ell$ is 
$$
\left<x_i\ell:i=1,\ldots, d\right>.
$$

\begin{proof}[Proof of Proposition \ref{ql_red}]
Suppose that for any rational linear form $l$ from $p(\Omega)$, the quadratic form 
$(\alpha Q+\beta L^2)|_{l=0}$ is rational for some $\alpha$ and $\beta$ with $(\alpha,\beta)\ne 0$.
This means that $F_l (\left<Q,L^2\right>)$
has a rational subspace of codimension $1$, or, equivalently, for some $\alpha,\beta\in\mathbb{R}$
such that $(\alpha,\beta)\ne 0$,
the space $\left<\alpha Q+\beta L^2,x_il:i=1,\ldots, d\right>$
is rational.

{\sc Case 1:} For any rational linear forms $l$ and $m$ such that $(l,m)\in\pi(\Omega)$, the space 
$\left<Q,L^2,x_il,x_jm:i,j=1,\ldots, d\right>$ is rational. 

We will use the following claim:

\vspace{0.1cm}
\noindent {\sc Claim.} {\it Let $S$ be not in $\left<Q,L^2\right>$. Then there are rational
linear forms $l$ and $m$ such that $(l,m)\in\pi(\Omega)$, and $S$ is not in
$$
\left<Q,L^2,x_il,x_jm:i,j=1,\ldots, d\right>.
$$
}
\begin{proof}
Suppose not. Then for any rational linear forms $l$ and $m$ such that $(l,m)\in\pi(\Omega)$,
we have $S=\alpha(l,m)Q+\beta(l,m)L^2$ on the subspace $\{l=0\}\cap\{m=0\}$
for some $\alpha(l,m),\beta(l,m)\in\mathbb{R}$.
First, we show that $\alpha(l,m)$ and $\beta(l,m)$ are independent of $(l,m)$
on a large set of pairs. Let $(l,m,n)\in\Omega$ for some rational linear forms $l$, $m$, $n$.
Because of the condition (3), $Q$ has rank at least $2$ on the space 
$$
V\stackrel{def}{=}\{l=0\}\cap\{m=0\}\cap\{n=0\}.
$$
Therefore, $Q|_V$ and $L^2|_V$ are linearly independent.
It follows that
\begin{equation} \label{ql_help1}
\alpha(l,m)=\alpha(m,n)=\alpha(n,l)\quad\hbox{and}\quad \beta(l,m)=\beta(m,n)=\beta(n,l).
\end{equation}
Fix $(l_0,m_0)\in\pi(\Omega)$ with rational linear forms $l_0$ and $m_0$.
Consider the set 
$$
\mathcal{O}=\{(l,m)\in\pi(\Omega):(m,m_0)\in \pi(\Omega)\}.
$$
This is a nonempty open set in $\pi(\Omega)$. 
By Lemma \ref{ql_red1}, there is a rational linear form $n$ such that 
$(l_0,m_0,n),(l,m,n),(m,m_0,n)\in\Omega$. Then using (\ref{ql_help1}), we obtain
$$
\alpha(l_0,m_0)=\alpha(m_0,n)=\alpha(m,n)=\alpha(l,m),
$$
and similarly, $\beta(l_0,m_0)=\beta(l,m)$.
Hence, the coefficients $\alpha(l,m)$ and $\beta(l,m)$ are constant for rational linear forms $l$ and $m$
such that $(l,m)\in\mathcal{O}$.
Thus, $S-\alpha Q-\beta L^2 = 0$ on the space $\{l=0\}\cap\{m=0\}$ for rational linear forms $l$ and $m$
such that $(l,m)\in\mathcal{O}$. It follows that $S-\alpha Q-\beta L^2 = 0$
on an open subset of $\mathbb{R}^d$. Hence, $S=\alpha Q+\beta L^2$.
This is a contradiction.
\end{proof}

Using the claim, we conclude that the space $\left<Q,L^2\right>$ is an intersection of spaces of
the form $\left<Q,L^2,x_il,x_jm:i,j=1,\ldots n\right>$ for some rational linear forms $l$ and $m$ such that
$(l,m)\in\mathcal{O}$. By the assumption of Case 1, each of these spaces is rational.
Hence, the space $\left<Q,L^2\right>$ is rational too. This gives a contradiction.

{\sc Case 2:} For some rational linear forms $l$ and $m$ such that 
$(l,m)\in\pi (\Omega )$, the space $\left<Q,L^2,x_il,x_jm:i,j=1,\ldots, d\right>$ is not rational. 

By the assumption, the space $\left<\alpha Q+\beta L^2,x_il:i=1,\ldots, d\right>$
is rational for some $\alpha,\beta\in\mathbb{R}$ such that $(\alpha,\beta)\ne 0$.
Then the space 
$$
\left<\alpha Q+\beta L^2,x_il,x_jm:i,j=1,\ldots, d\right>
$$
is rational too.
It follows from the assumption of the Case 2 that the pair $(\alpha,\beta)$ with this property
is uniquely defined up to a scalar multiple. Let $\tilde{Q}=\alpha Q+\beta L^2$.
We show that $\tilde{Q}$ is proportional to a rational form. The pair $(\alpha,\beta)$ can be
chosen such that the form $\tilde{Q}|_{l=0}$ is rational.

Because of the uniqueness of $(\alpha,\beta)$, $\tilde{Q}|_{k=0}$ is proportional to a rational form
for any rational linear form $k\in p(\Omega)\cap\left<l,m\right>$.
Note that the later set is not empty because $l$ and $m$ is in it.
Let $\{ e_i: i=1,\ldots, d-2\}$ be a rational basis of $\{l=0\}\cap\{m=0\}$. We complete
it to a rational basis of $\mathbb{R}^d$ by vector $e_{d-1}$ and $e_d$ such that $l(e_{d-1})=0$.
Let 
$$
U_t=\left<e_i,e_{d-1}+te_d:i=1,\ldots, d-2\right>\;\;\hbox{for}\;\;t\in\mathbb{Q}.
$$
Since $U_0=\{l=0\}$ and $l\in p(\Omega )$, a rational form $k_t$ that defines
$U_t$ is in $p(\Omega )$ for sufficiently small $t$. Also $k_t\in\left<l,m\right>$ because
$$
\{l=0\}\cap\{m=0\}\subseteq U_t.
$$
Therefore, $\alpha_t \tilde{Q}|_{U_t}$ is rational for some $\alpha_t\ne 0$ when $t$ is rational and
sufficiently small. Since $(l,m)\in \pi(\Omega)$, $\tilde{Q}\ne 0$ on the space $\{l=0\}\cap\{m=0\}$.
Take a rational vector $x\in\{l=0\}\cap\{m=0\}$ such that $\tilde{Q}(x)\ne 0$.
Then $\tilde{Q}(x)\in\mathbb{Q}$, and $\alpha_t \tilde{Q}(x)\in\mathbb{Q}$. Therefore, $\alpha_t\in\mathbb{Q}$,
and $\tilde{Q}|_{U_t}$ is rational for sufficiently small $t$.
In particular, when $t=0$, $\tilde{Q}(e_i,e_j)\in\mathbb{Q}$
for $i,j<d$. Also 
$$
\tilde{Q}(e_i,e_{d-1}+te_d)=\tilde{Q}(e_i,e_{d-1})+t\tilde{Q}(e_i,e_d)\in\mathbb{Q},
$$
so that $\tilde{Q}(e_i,e_d)\in\mathbb{Q}$ for $i<d$. Finally,
$$
\tilde{Q}(e_{d-1}+te_d)=\tilde{Q}(e_{d-1})+2t\tilde{Q}(e_{d-1},e_d)+t^2\tilde{Q}(e_d)\in\mathbb{Q}.
$$
Thus, $\tilde{Q}(e_d)\in\mathbb{Q}$. This shows that $\tilde{Q}$ is rational which is a contradiction. 
\end{proof}

\section{Pairs of type (I)}

In this section we consider a pair $(Q,L)$ of type (I) and dimension $4$. By Proposition \ref{ql_can}, the pair 
$(Q,L)$ is equivalent to either 
$$
(x_1^2+x_2^2-x_3^2-x_4^2, x_4)\;\;\hbox{or}\;\; (x_1^2-x_2^2-x_3^2-x_4^2, x_4).
$$
We consider the first case. The other case can be done the same way.
Let 
\begin{equation} \label{ql_canp}
Q_0=x_1^2+x_2^2-x_3^2-x_4^2\quad\hbox{and}\quad L_0=x_4.
\end{equation}


Let
\begin{equation} \label{ql_st}
H=\left(
\begin{tabular}{c|c}
$\hbox{SO(2,1)}$ & $0$\\
\hline
$0$              & $1$
\end{tabular}
\right) \subseteq \hbox{SL}(4,\mathbb{C}).
\end{equation}
Note that the group $H$ leaves the pair $(Q_0, L_0)$ invariant.

First, we collect some simple facts.

\begin{pro} \label{ql_subgr}
Let $F$ be a connected semisimple algebraic subgroup of $\hbox{SL}(4,\mathbb{C})$
which acts irreducibly on $\mathbb{C}^4$ such that
$$
H\subsetneq F\subsetneq \hbox{SL}(4,\mathbb{C}).
$$
Then $G=\hbox{SO}(S,\mathbb{C})$ for some nondegenerate quadratic form $S$.
\end{pro}

\begin{proof}
Let $\mathfrak f$ be the Lie algebra of $F$.
For convenience of the reader, let us reproduce in Table 1
the list of complex semisimple Lie algebras of dimension
up to $14$ and dimensions of their faithful irreducible representations (see \cite{OV}).
Note that an irreducible representation of a semisimple algebra is a tensor product
of irreducible representations of simple factors.
\begin{table}[ht]
\begin{center}
\begin{tabular}{|c|c|l|}
\hline
$\mathfrak g$ & $\dim \mathfrak g$ & $\dim \rho$\\
\hline
$A_1$             & $3$  & $2,3,4,\ldots$\\
$A_1+A_1$         & $6$  & $4,\ldots$\\
$A_1+A_1+A_1$     & $9$  & $8,\ldots$\\
$A_1+A_1+A_1+A_1$ & $12$ & $16,\ldots$\\
$A_2$             & $8$  & $3,8,\ldots$\\
$A_2+A_1$         & $11$ & $6,\ldots$\\
$A_2+A_1+A_1$     & $14$ & $12,\ldots$\\
$B_2=C_2$         & $10$ & $4,5,\ldots$\\
$B_2+A_1$         & $13$ & $8,\ldots$\\
$G_2$             & $14$ & $7,\ldots$\\
\hline
\end{tabular}
\end{center}
\caption{Complex irreducible representations}
\end{table}
One can see that $\mathfrak f$ should be of type $A_1+A_1$ or $C_2$.
The $4$-dimensional irreducible representation of the algebra of type $C_2$ is simplectic.
On the other hand, $H\subset F$ acts irreducibly on a $3$-dimensional subspace,
and leaves invariant a nonzero symmetric form on this subspace.
This shows that $F$ cannot leave a nonzero symplectic form invariant.
Thus, $\mathfrak f$ is of type $A_1+A_1$.
Since $\mathfrak f$ has a unique irreducible $4$-dimensional representation,
this representation is equivalent to the representation of $\mathfrak{so}(4,\mathbb{C})$
on $\mathbb{C}^4$.
This means that $\mathfrak f=g\mathfrak{so}(4,\mathbb{C})g^{-1}$ for
some $g\in\hbox{SL}(4,\mathbb{C})$. Equivalently, $\mathfrak f=\mathfrak{so}(S,\mathbb{C})$
for some nondegenerate quadratic form $S$. 
Since $F$ is connected, $F=\hbox{SO}(S,\mathbb{C})$. 
\end{proof}


\begin{lem} \label{ql_ss}
Let $G\subseteq\hbox{GL}(d,\mathbb{C})$ be an algebraic group that acts irreducibly on $\mathbb{C}^d$.
Then its unipotent radical is trivial.
\end{lem}

\begin{proof} Let $U$ be the unipotent radical of $G$, and $V=\mathbb{C}^d$.
Consider the space 
$$
V^U=\{x\in V:Ux=x\}.
$$
Since $U$ is unipotent, $V^U\ne 0$. Since $U$ is normal in $G$, $V^U$ is $G$-invariant.
Thus, $V^U=V$, and $U=1$.
\end{proof}

\begin{lem} \label{ql_inv}
Let $Q$ be a quadratic form on $\mathbb{C}^4$ which is invariant
under the action of $H$ ($H$ is defined in (\ref{ql_st})). Then $Q=\alpha Q_0+\beta L_0^2$
for some $\alpha,\beta\in\mathbb{C}$, where $Q_0$ and $L_0$ are defined in (\ref{ql_canp}).
\end{lem}

\begin{proof}
Let us write $Q=Q(x_1,x_2,x_3,0)+L(x_1,x_2,x_3)x_4+\beta x_4^2$
for some linear form $L$ and $\beta\in\mathbb{C}$. Then $Q(\cdot,\cdot,\cdot,0)$ and $L$ is
$\hbox{SO}(2,1)$-invariant. It follows that $Q(\cdot,\cdot,\cdot,0)=\alpha Q_0(\cdot,\cdot,\cdot,0)$
for some $\alpha\in\mathbb{C}$ and $L=0$.
\end{proof}

The following proposition is the main result of this section and a partial case
of Theorem \ref{ql_m}.

\begin{pro} \label{ql_I}
Let $Q$ be a nondegenerate quadratic form, and $L$ be a linear form
such that $(Q(x),L(x))=(\lambda Q_0(gx), \mu L_0(gx))$ for some 
$g\in \hbox{SL}(4,\mathbb{R})$ and $\lambda,\mu\in\mathbb{R}-\{0\}$, where $Q_0$ and $L_0$ are defined in (\ref{ql_canp}).
Suppose that for every $(\alpha,\beta)\ne (0,0)$, with $\alpha,\beta\in\mathbb{R}$,
the linear combination $\alpha Q+\beta L^2$ is not rational.

Then $\{(Q(x),L(x)):x\in\mathcal{P}(\mathbb{Z}^4)\}$ is dense in $\mathbb{R}^2$.
\end{pro}

\begin{proof}
Without loss of generality, $\lambda=\mu=1$.
Let $G=\hbox{SL}(4,\mathbb{C})$
and $H$ be as in (\ref{ql_st}). We study the action of $H_\mathbb{R}$ on $G_\mathbb{R}/G_\mathbb{Z}$.
By Ratner's Theorem \cite{r91},
the closure $\overline{H_\mathbb{R}gG_\mathbb{Z}}\subseteq G_\mathbb{R}/G_\mathbb{Z}$ is a homogeneous space.
Moreover by the result of Shah \cite[Proposition 3.2]{sh91}, $\overline{H_\mathbb{R}gG_\mathbb{Z}}=gF^o_\mathbb{R}G_\mathbb{Z}$
where $F$ is the smallest algebraic $\mathbb{Q}$-subgroup containing $g^{-1}Hg$,
and the radical of $F$ is unipotent.
Here $F^o_\mathbb{R}$ denotes the connected component of $F_\mathbb{R}$ with respect to 
the Euclidean topology. Since $H$ is connected as an algebraic group, $F$ is connected too.

First, we consider the case when $F$ acts irreducibly on $\mathbb{C}^4$.
By Lemma \ref{ql_ss}, $F$ is semisimple, and by Proposition \ref{ql_subgr},
$F$ is one of the subgroups $g^{-1}Hg$, $G$, $\hbox{SO}(S,\mathbb{C})$ for some nondegenerate quadratic form $S$.
Since $F$ acts irreducibly on $\mathbb{C}^4$, $F\ne g^{-1}Hg$.
Suppose that $F=\hbox{SO}(S)$. Since $F$ is defined
over $\mathbb{Q}$, the quadratic form $S$ is proportional to a rational form
(see, for example, \cite{bo95}).
By Lemma \ref{ql_inv}, $S=\alpha Q+\beta L^2$ for some
$\alpha,\beta\in\mathbb{C}$. This is a contradiction. Thus, $F=G$.
We conclude that $\overline{H_\mathbb{R}gG_\mathbb{Z}}=G_\mathbb{R}/G_\mathbb{Z}$.
Then
\begin{eqnarray*}
&&\overline{\{(Q(x),L(x)):x\in\mathcal{P}(\mathbb{Z}^4)\}}\\
&=&\overline{\{(Q_0(hg\gamma e_1),L_0(hg\gamma e_1)):h\in H_\mathbb{R},\gamma\in G_\mathbb{Z}\}}\\
&=&\{(Q_0(ae_1),L_0(ae_1):a\in G_\mathbb{R}\}=\{(Q_0(x),L_0(x):x\in \mathbb{R}^4\}=\mathbb{R}^2.
\end{eqnarray*}
Here $e_1=(1,0,\ldots,0)$, and the last equality holds by Lemma \ref{ql_im}.

Now we assume that the action of $F$ on $\mathbb{C}^4$ is not irreducible.
Then the action of $F$ on the dual space $\mathcal{L}$ is reducible too.
The dual action of $H$ has irreducible components $\left<x_1,x_2,x_3\right>$
and $\left<x_4\right>$. Moreover, for any linear form $l\notin \left<x_4\right>$,
the smallest $H$-invariant subspace containing $l$ should contain
$\left<x_1,x_2,x_3\right>$. Therefore, $\left<x_1,x_2,x_3\right>$ and $\left<x_4\right>$ are the only nontrivial $H$-invariant subspaces.
For $i=1,\ldots,4$, define linear forms $\ell_i(x)=(gx)_i$. Then
$$
Q=\ell_1^2+\ell_2^2-\ell_3^2-\ell_4^2\quad\hbox{and}\quad L=\ell_4.
$$
It follows that $\mathcal{L}_1\stackrel{def}{=}\left<\ell_1,\ell_2,\ell_3\right>$ and $\mathcal{L}_2\stackrel{def}{=}\left<\ell_4\right>$
are the only candidates for nontrivial $F$-invariant subspaces.
Let $F$ be semisimple. Then its action is completely reducible.
Therefore, $\mathcal{L}_1$ and $\mathcal{L}_2$ are $F$-invariant, and since $F$ has no nontrivial characters,
the space of $F$-fixed vectors
$\mathcal{L}^F=\mathcal{L}_2=\left<L\right>$. Using that $F$ is a $\mathbb{Q}$-group, we conclude that $\left<L\right>$ is a $\mathbb{Q}$-subspace.
This means that $L$ is proportional to a rational form, which is a contradiction.
Therefore, $F$ is not semisimple, and its unipotent radical $U\ne 1$. Note that $U$ is defined over $\mathbb{Q}$.
The subspace of $U$-fixed vectors $\mathcal{L}^U$ is $F$-invariant, so that either
$\mathcal{L}^U=\mathcal{L}_1$ or $\mathcal{L}^U=\mathcal{L}_2=\left<L\right>$. If the later holds,
$L$ is proportional to a rational form by the same argument as above.
Thus, we may assume that $\mathcal{L}^U=\mathcal{L}_1$. In particular, $\mathcal{L}_1$ is a $\mathbb{Q}$-subspace.
Since $F$ is defined over $\mathbb{Q}$, it has a Levi subgroup $F_0$ which is
defined over $\mathbb{Q}$. It follows from the description of $F$-invariant subspaces
that $F_0$ acts irreducibly on $\mathcal{L}^U$. There exists an $F_0$-invariant complement for $\mathcal{L}^U$
generated by a linear form $\ell_0$. We have that $\mathcal{L}^{F_0}=\left<\ell_0\right>$. In particular, 
$\left<\ell_0\right>$ is a $\mathbb{Q}$-subspace, and the linear form $\ell_0$ can be chosen to rational.
By Malcev's theorem, $g^{-1}Hg\subseteq u^{-1}F_0u$ for some $u\in U$.
The forms $L(u^{-1}\cdot x)$ and $\ell_0(x)$ are both fixed by $F_0$,
so that $L(u^{-1}\cdot x)=\alpha\cdot \ell_0(x)$ for some $\alpha\in\mathbb{R}-\{0\}$.

Suppose that $g^{-1}Hg=u^{-1}F_0u$. Then $F_0$ stabilizes the pair
$$
\left(Q(u^{-1}\cdot x),L(u^{-1}\cdot x)\right)=\left(\ell_1^2+\ell_2^2-\ell_3^2-\alpha^2 \ell_0^2,\alpha \ell_0\right)
$$
Since $F_0$ is a $\mathbb{Q}$-group, the space of quadratic forms that are fixed by $F_0$ is
a $\mathbb{Q}$-space. By Lemma \ref{ql_inv}, this space is spanned by
$\ell_1^2+\ell_2^2-\ell_3^2-\alpha^2 \ell_0^2$ and $\ell_0^2$. Therefore, for some $\alpha,\beta\in\mathbb{R}$,
the form 
$$
\tilde{Q}\stackrel{def}{=}\alpha(\ell_1^2+\ell_2^2-\ell_3^2)+\beta \ell_0^2
$$
is rational. Since the space
$\mathcal{L}_1=\left<\ell_1,\ell_2,\ell_3\right>$ is rational, there exists a rational vector $x_0\ne 0$
such that $\ell_1(x_0)=\ell_2(x_0)=\ell_3(x_0)=0$. Clearly, $\ell_0(x_0)\ne 0$.
Then $\tilde{Q}(x_0)=\beta\cdot \ell_0(x_0)^2\in\mathbb{Q}$. Therefore, $\beta\in\mathbb{Q}$,
and the form $Q+L^2=\ell_1^2+\ell_2^2-\ell_3^2$ is proportional to a rational form,
which is a contradiction.

Now we may assume that $g^{-1}Hg\subsetneq u^{-1}F_0u$. Then
$$
F_0=ug^{-1}\left(
\begin{tabular}{c|c}
$\hbox{SL}(3)$ & $0$\\
\hline
$0$              & $1$
\end{tabular}
\right)gu^{-1}.
$$
Also
$$
F_1\stackrel{def}{=}g^{-1}\left(
\begin{tabular}{c|c}
$\hbox{SL}(3)$ & $0$\\
\hline
$0$              & $1$
\end{tabular}
\right)g\subseteq F.
$$
An element $u\in U$ acts on $\mathcal{L}$ as follows: $u\cdot \ell=\ell$ for $\ell\in \mathcal{L}_1$, and
$u\cdot \ell_4=\ell_4+\ell_u$ for some linear form $\ell_u\in \mathcal{L}_1$.
Then the adjoint action of $F_1$ on $U$ corresponds to the usual action of $F_1$
on the space spanned by linear forms $\ell_u$, $u\in U$. It follows that for every $\ell\in \mathcal{L}_1$, there exists $u\in U$ such that
$u\cdot \ell_4=\ell_4+\ell$. 

Let $\{e_i:i=1,\ldots,4\}$ be the standard basis of $\mathbb{C}^4$.
Fix $(a,b)\in\mathbb{R}^2$. Take $x\in\mathcal{P}(\mathbb{Z}^4)-\left<g^{-1}e_4\right>$.
By the previous remark, there exists $u\in U_\mathbb{R}$ such that
$$
L(ux)=(u\cdot \ell_4)(x)=b.
$$
Write $ux=x_1+x_2$ for $x_1\in g^{-1}\left<e_1,e_2,e_3\right>$ and
$x_2\in g^{-1}\left<e_4\right>$. 
Since $Q$ is indefinite on the subspace $g^{-1}\left<e_1,e_2,e_3\right>$,
there exists $f\in F_1$ such that $Q(fx_1)=a-Q(x_2)$.
Then 
$$
Q(fux)=Q(fx_1)+Q(x_2)=a.
$$
Since 
$$
\overline{g^{-1}H_\mathbb{R}gG_\mathbb{Z}}=F^o_\mathbb{R}G_\mathbb{Z}=F_{1\mathbb{R}}U_\mathbb{R}G_\mathbb{Z},
$$
there exist $h_n\in g^{-1}H_\mathbb{R}g$ and $\gamma_n\in G_\mathbb{Z}$ such that
$h_n\gamma_n\rightarrow fu$ as $n\rightarrow\infty$.
Finally,
\begin{eqnarray*}
(Q(\gamma_nx),L(\gamma_n x))&=&(Q(h_n\gamma_nx),L(h_n\gamma_n x))\\
&\rightarrow& (Q(fux),L(fux))=(a,b).
\end{eqnarray*}
The proposition is proved.
\end{proof}


\begin{remark}
Proposition \ref{ql_I} combined with Proposition \ref{ql_red} implies Theorem \ref{ql_m} for dimension $d\ge 5$.
\end{remark}

\section{Pairs of type (II)}

Now we prove Theorem \ref{ql_m} for pairs of type (II) of dimension $4$.
This will finish the proof of Theorem \ref{ql_m}. Let
\begin{equation} \label{ql_canII}
Q_0=x_1x_4-x_2x_3\quad\hbox{and}\quad L_0=x_4.
\end{equation}
Note that we use different (cf. Proposition \ref{ql_can}) canonical form to
simplify calculations.

Let
\begin{equation} \label{ql_stII}
H=\left\{\left(
\begin{tabular}{cccc}
$1$ & $a$ & $b$ & $ab$\\
$0$ & $1$ & $0$ & $b$\\
$0$ & $0$ & $1$ & $a$\\
$0$ & $0$ & $0$ & $1$
\end{tabular}
\right): a,b\in\mathbb{C}\right\} \subseteq \hbox{SL}(4,\mathbb{C}).
\end{equation}
It is easy to check that $Q_0(hx)=Q_0(x)$ and $L_0(hx)=L_0(x)$ for $h\in H$.

\begin{lem} \label{ql_invII}
The only nontrivial $H$-invariant subspaces of the dual action of $H$ 
(i.e the action on the space $\mathcal{ L}$ of linear forms) are
\begin{itemize}
\item $\left<x_4\right>$.
\item $\left<\alpha x_2+\beta x_3,x_4\right>$ for some $\alpha,\beta\in \mathbb{C}$.
\item $\left<x_2, x_3,x_4\right>$.
\end{itemize}
\end{lem}

\begin{proof}
One can check that $x_4$ is the only fixed vector of $H$ (up to a scalar multiple).
Since the action of $H$ is unipotent, every nontrivial $H$-invariant subspace $\mathcal{ V}$
contains a nonzero vector fixed by $H$. It follows that $\left<x_4\right>\subseteq \mathcal{ V}$.
Consider a factor space $\mathcal{ L}/\left<x_4\right>$.
The subspace $\left<x_2, x_3,x_4\right>/\left<x_4\right>$ consists of
$H$-fixed vectors,  and for any $v\notin\left<x_2, x_3,x_4\right>$,
one has $\left<Hv,x_4\right>=\mathcal{L}$. The conclusion follows.
\end{proof}

\begin{lem} \label{ql_prod}

\begin{enumerate}
\item[1.] Let $G$ be a simple group, and $\bar{G}=G\times G$. Denote by $\pi_1$ and $\pi_2$
the projection maps. Let $H\lneq \bar{G}$ be such that $\pi_i(H)=G$ for $i=1,2$.
Then $H=\{(\alpha(g),g):g\in G\}$ for some automorphism $\alpha$ of $G$.

\item[2.] Let $\mathfrak{g}=\mathfrak{sl}(2,\mathbb{C})$, and $\bar{\mathfrak{g}}=\mathfrak{g}\oplus\mathfrak{g}$.
Denote by $\pi_1$ and $\pi_2$
the projection maps. Let $\mathfrak{h}\lneq {\bar{\mathfrak{g}}}$ be such that
$\pi_i(\mathfrak{h})=\mathfrak{g}$ for $i=1,2$.
Then $\mathfrak{h}=\{(\hbox{Ad}(g)x,x):x\in \mathfrak{g}\}$ for some $g\in \hbox{SL}(2,\mathbb{C})$.
\end{enumerate}
\end{lem}

\begin{proof}
Let 
$$
S_1=\{g\in G:(g,e)\in H\}\;\;\hbox{and}\;\; S_2=\{g\in G:(e,g)\in H\}.
$$
It is easy to check that $S_1$ and $S_2$ are normal subgroups of $G$, that is,
$g^{-1}S_ig\subseteq S_i$ for every $g\in G$, $i=1,2$.
Since $G$ is simple, $S_i$ is either $\{e\}$ or $G$. 

Suppose that $S_1=G$. We show that $H=\bar{G}$. Let $(g,h)\in \bar{G}$. For some $g_1\in G$,
$(g_1,h)\in H$. Then 
$$
(g,h)=(g_1,h)(g_1^{-1} g,e)\in H.
$$
Hence, $H=\bar{G}$.
Similarly, $H=\bar{G}$ if $S_2=G$.

Suppose that $S_1=S_2=\{e\}$. For any $g\in G$, there exists a unique element $\alpha (g)$
such that $(\alpha (g),g)\in H$. Because of the uniqueness, $\alpha$ is a homomorphism.
It is surjective because $\pi_1(H)=G$ and injective because $S_2=\{e\}$.
This proves the first part of the lemma.

It is straightforward to rewrite this argument for simple Lie algebras. It is known
that any automorphism of $\mathfrak{sl}(2,\mathbb{C})$ is inner.
Thus, the second part of the lemma follows.
\end{proof}


The following proposition finishes the proof of Theorem \ref{ql_m}.
Its proof is similar to the proof of Proposition \ref{ql_I}.
\begin{pro} \label{ql_II}
Let $Q$ be a nondegenerate quadratic form, and $L$ be a linear form
such that $(Q(x),L(x))=(\lambda Q_0(gx), \mu L_0(gx))$ for some 
$g\in \hbox{SL}(4,\mathbb{R})$ and $\lambda,\mu\in\mathbb{R}-\{0\}$, where $Q_0$ and $L_0$ are defined in (\ref{ql_canII}).
Suppose that 
for every $(\alpha,\beta)\ne (0,0)$, with $\alpha,\beta\in\mathbb{R}$,
the linear combination $\alpha Q+\beta L^2$ is not rational.

Then $\{(Q(x),L(x)):x\in\mathcal{P}(\mathbb{Z}^4)\}$ is dense in $\mathbb{R}^2$.
\end{pro}

\begin{proof}
We may assume that $\lambda=\mu=1$.
Let $G=\hbox{SL}(4,\mathbb{C})$,
and $H$ be as in (\ref{ql_stII}). Consider the action of $H_\mathbb{R}$ on $G_\mathbb{R}/G_\mathbb{Z}$.
As in the proof of Proposition \ref{ql_I}, we know that
$\overline{H_\mathbb{R}gG_\mathbb{Z}}=gF^o_\mathbb{R}G_\mathbb{Z}$
where $F$ is the smallest algebraic $\mathbb{Q}$-subgroup containing $g^{-1}Hg$.
The group $F$ is connected as an algebraic group, and
its radical is unipotent. There is Levi decomposition
$F=F_0U$ where $F_0$ is a connected (as an algebraic group) semisimple $\mathbb{Q}$-subgroup, and
$U$ is the unipotent radical of $F$. Note that $U$ is defined over $\mathbb{Q}$.

To prove the proposition, it is enough to show that
\begin{equation} \label{ql_eqII}
\overline{\{(Q(fz),L(fz)):f\in F^o_\mathbb{R}, z\in \mathcal{P}(\mathbb{Z}^4) \}}=\mathbb{R}^2.
\end{equation}

Consider the action of $F$ on the space of linear forms $\mathcal{ L}$.
Let $\mathcal{ V}$ be a nontrivial $F$-invariant subspace. 
Since $H$ is unipotent, $\mathcal{V}$ contains a nonzero vector fixed by $g^{-1}Hg$.
By Lemma \ref{ql_invII}, the only vector fixed by $g^{-1}Hg$ is $L$
(up to a scalar multiple). Thus, $\left< L\right>\subseteq \mathcal{ V}$.
It follows that there is a unique $F$-irreducible subspace $\mathcal{ V}\subseteq \mathcal{ L}$.
Namely, it is the intersection of all nontrivial $F$-invariant subspaces.
This subspace is contained in every $F$-invariant subspace.
Also $\mathcal{ V}$ is defined over $\mathbb{Q}$. Indeed, for any 
$\sigma\in\hbox{Gal}(\mathbb{C}/\mathbb{Q})$, $\mathcal{ V}^\sigma$ is 
$F(\mathbb{Q})$-invariant, and $F(\mathbb{Q})$ is Zariski dense in $F$
because $F$ is connected. Therefore, $\mathcal{ V}^\sigma$ is $F$-invariant,
and $\mathcal{ V}\subseteq \mathcal{ V}^\sigma$. Comparing dimensions, we conclude that $\mathcal{ V}=\mathcal{ V}^\sigma$.
This shows that $\mathcal{V}$ is defined over $\mathbb{Q}$.

Let $\mathcal{ L}^U$ be the subspace of vectors fixed by $U$.
Since $U$ is normal in $F$, the space $\mathcal{ L}^U$ is $F$-invariant,
so that $\mathcal{ V}\subseteq \mathcal{ L}^U$. It follows that $F_0$ act irreducibly on $\mathcal{ V}$.
Since $F_0$ is semisimple, the $F$-action on $\mathcal{ L}^U$ is completely reducible.
Suppose that $\mathcal{ V}\subsetneq \mathcal{ L}^U$. Then $\mathcal{ L}^U=\mathcal{ V}\oplus \mathcal{ W}$
for some $F$-invariant subspace $\mathcal{ W}$. However, this contradicts
the description of $H$-invariant subspaces in Lemma \ref{ql_invII}.
Thus, $\mathcal{ V}=\mathcal{ L}^U$.

We can write
$$
\mathcal{ L}_\mathbb{R}=\mathcal{ V}_\mathbb{R}\oplus \mathcal{ W}_\mathbb{R}
$$
for some real $(F_0)_\mathbb{R}$-invariant subspace $\mathcal{ W}_\mathbb{R}$.
Then $\mathcal{ L}=\mathcal{ V}\oplus \mathcal{ W}$ where 
$\mathcal{ W}=\mathcal{ W}_\mathbb{R}\otimes \mathbb{C}$ is $F_0$-invariant,
because $(F_0)_\mathbb{R}$ is Zariski dense in $F_0$.
Note that $\mathcal{ W}$ is defined over $\mathbb{R}$.
Let 
$$
V=\{v\in\mathbb{C}^4: \ell(v)=0\;\;\hbox{for}\;\ell\in\mathcal{ V} \}
$$
and 
$$
W=\{v\in\mathbb{C}^4: \ell(v)=0\;\;\hbox{for}\;\ell\in\mathcal{ W} \}.
$$
Clearly, $\mathbb{C}^4=V\oplus W$, $V$ is $F$-invariant, and $W$ is
$F_0$-invariant. Moreover, $W$ is $F_0$-irreducible because $W\simeq\mathbb{C}^4/V$
as $F_0$-modules, and any nontrivial $F_0$-invariant subspace
of $\mathbb{C}^4/V$  would give by duality a nontrivial $F_0$-invariant subspace in $\mathcal{ V}$.
The space $(\mathbb{C}^4/V)^U$ is nonzero and $F_0$-invariant.
It follows that $(\mathbb{C}^4/V)^U=\mathbb{C}^4/V$, so that $U$ acts trivially on $\mathbb{C}^4/V$.

Let $\ell_i(x)=(gx)_i$ for $i=1,\ldots,4$. Then $Q=\ell_1\ell_4-\ell_2\ell_3$ and $L=\ell_4$.

Consider several cases:

\underline{{\sc Case 1:} $\dim \mathcal{ V}=1$.} Then $\mathcal{ V}=\left< L\right>$,
and since $\mathcal{V}$ is defined over $\mathbb{Q}$, $L$ is a multiple of a rational form.
This is a contradiction.

\underline{{\sc Case 2:} $\dim \mathcal{ V}=2$.}
Then $\dim V=\dim W=2$.
Since $F_0$ is semisimple, the action of $F_0$ on  $V$
is either trivial or irreducible.

Denote 
$$
A=\{g\in\hbox{SL}(4,\mathbb{C}):g|_V=id,g|_W=id+a\;\hbox{for}\; a\in\hbox{End}(W,V)\}.
$$

Suppose that $F_0$ acts irreducibly on $V$.
Then $U$ acts trivially on $V$ because otherwise $V^U\ne 0$ is a nontrivial
$F_0$-invariant subspace. Since it was shown above that $U$ acts trivially on $\mathbb{C}^4/V$ too,
this implies that $U\subseteq A$.
The Lie algebra $\mathfrak{f}_0\subseteq\mathfrak{sl}(V)\times\mathfrak{sl}(W)$ of $F_0$
satisfies the conditions of Lemma \ref{ql_prod} unless $\mathfrak{f}_0=\mathfrak{sl}(V)\times\mathfrak{sl}(W)$.

First, we consider the case when the last equality holds. Then $A$ is an irreducible $F_0$-module.
Suppose that $U=1$. Then $\mathcal{V}$ and $\mathcal{W}$ are $F$-invariant subspaces
such that $\mathcal{V}+\mathcal{W}=\mathcal{L}$.
Since every two $2$-dimensional subspaces in Lemma \ref{ql_invII} are contained in the
unique $3$-dimensional subspace, this gives a contradiction. Thus, $U\ne 1$, and since it is a submodule
of the irreducible $F_0$-module $A$, we conclude that $U=A$. In particular, for any $x\notin V_\mathbb{R}$ and any $v\in V_\mathbb{R}$,
there exists $u\in U_\mathbb{R}$ such that $ux=x+v$. We will use this fact latter.

Now we assume that $\mathfrak{f}_0$ is a proper subalgebra of $\mathfrak{sl}(V)\times\mathfrak{sl}(W)$.
Then by Lemma \ref{ql_prod},
$$
\mathfrak{f}_0=\{(\phi^{-1}x\phi,x):x\in\mathfrak{sl}(W)\}
$$
for some isomorphism $\phi:V\rightarrow W$.
Since $\mathfrak f_0$ is an $\mathbb{R}$-subalgebra, $\phi$ can be taken
to be an $\mathbb{R}$-map.
Let $\mathfrak a\simeq\{a\in\hbox{End}(W,V)\}$ be the Lie algebra of $A$.
The adjoint action of $\mathfrak{f}_0$ on $\mathfrak a$
is isomorphic to adjoint action of $\mathfrak{sl}(W)$ on
$\mathfrak{gl}(W)$. The isomorphism is
$a\rightarrow \phi\circ a$. It follows that $\mathfrak a$ has the
only nontrivial $\mathfrak f_0$-modules: $\{\lambda \phi^{-1}: \lambda\in\mathbb{C}\}$
and $\{ \phi^{-1}\circ b :b\in \mathfrak{sl}(W)\}$.

Let $U$ correspond to $\{\lambda \phi^{-1}: \lambda\in\mathbb{C}\}$.
Take a basis $\{v_1,v_2\}$ of $V$. Then $\{v_1,v_2,\phi(v_1),\phi(v_2)\}$ is a basis
of $\mathbb{C}^4$, and with respect to this basis,
$$
F=\left\{\left(
\begin{tabular}{cc}
$g$ & $\lambda g$\\
$0$              & $g$
\end{tabular}
\right): g\in\hbox{SL}(2), \lambda\in \mathbb{C}\right\}.
$$
Let $S$ be an $F$-invariant quadratic form. Then the subspace $V$ and $W$ are totally isotropic
with respect to $S$. Thus, the matrix of $S$ is of the form
$\left(
\begin{tabular}{cc}
$0$ & $X$\\
$X$ & $0$
\end{tabular}
\right)$
for some matrix $X$ such that ${}^tgXg=X$ for all $g\in\hbox{SL}(2)$.
It follows that $X=\left(
\begin{tabular}{cc}
$0$ & $-u$\\
$u$ & $0$
\end{tabular}
\right)$
for some $u\in\mathbb{C}$, and the quadratic form $S$ is unique up to a scalar multiple.
Since $F$ is a $\mathbb{Q}$-group, the space of $F$-invariant quadratic forms is defined
over $\mathbb{Q}$. Therefore, $S$ is a multiple of a rational form, and 
$g^{-1}Hg$ is contained in the $\mathbb{Q}$-group $\hbox{SO}(S)$. Now we can argue as in the Case 4
when $F$ is of type $A_1+A_1$ (see below) to get a contradiction.
Hence, we may assume that $U$ corresponds
to $\{ \phi^{-1}\circ b :b\in \mathfrak{sl}(W)\}$.
Then for any $x\notin V_\mathbb{R}$ and any $v\in V_\mathbb{R}$,
there exists $u\in U_\mathbb{R}$ such that $ux=x+v$. We will use this fact latter on.

Suppose that $F_0$ acts trivially on $V$. If $U$ acts trivially on $V$, $F$
has linearly independent $F$-invariant vectors. They correspond to distinct
$3$-dimensional $F$-invariant subspaces in $\mathcal{ L}$. This contradicts
Lemma \ref{ql_invII}. Thus, $U$ acts nontrivially on $V$. Let $v_1\in V$
be a $U$-fixed vector in $V$. Take $v_2$ such that $V=\left<v_1,v_2\right>$.
Let 
$$
B=\{g\in\hbox{SL}(4,\mathbb{C}):gv_1=v_1,gv_2=v_2+tv_1,g|_W=id\;\hbox{for}\;t\in\mathbb{C}\}.
$$
Clearly, $U\subseteq AB$ and $U\nsubseteq A$.
If $U\subseteq B$, then the subspaces $V$ and $W$ would be $F$-invariant, which
contradicts Lemma \ref{ql_invII}. Thus, $U\nsubseteq B$.
Let $u=ab\in U$ for $a\in A-\{1\}$ and $b\in B-\{1\}$. Note that $A$ is abelian, and $F$ normalizes $A$.
Thus,
$$
u^{-1}(U\cap A)u=b^{-1}(U\cap A)b\subseteq U\cap A.
$$
Since the action of $B$ is algebraic, it follows that $B$ normalizes $U$. Let
$$
a=\left(
\begin{tabular}{cc}
$1$ & $\tilde{a}$\\
$0$              & $1$
\end{tabular}
\right)
$$
for some $\tilde{a}\in\hbox{End}(W,V)-\{0\}$ with respect to the decomposition $V\oplus W$.
Since $F_0$ acts irreducibly on $W$, the
restriction map is $F_0\rightarrow \hbox{SL}(W)$ is surjective.
For $g\in \hbox{SL}(W)$, take $f\in F_0$ such that $f=id\oplus g$ with respect to the decomposition
$V\oplus W$. Then
$$
f^{-1}uf=\left(
\begin{tabular}{cc}
$1$ & $\tilde{a}g$\\
$0$              & $1$
\end{tabular}
\right)\cdot b\in U.
$$
Thus,
$$
\left(
\begin{tabular}{cc}
$1$ & $\tilde{a}g-\tilde{a}$\\
$0$              & $1$
\end{tabular}
\right)\in U
$$
for every $g\in\hbox{SL}(W)$. In particular, it follows that $U\cap A\ne 1$.
The only nontrivial $BF_0$-submodule of $A$ is 
$$
A_0\stackrel{def}{=}
\left(
\begin{tabular}{c|c}
$1$ & $\hbox{End}(W,\left<v_1\right>)$\\
\hline
$0$              & $1$
\end{tabular}
\right).
$$
Suppose that $U\subseteq A_0B$. Then the linear form corresponding to the
projection to $\left<v_2\right>$ is fixed by $F$. This contradicts Lemma \ref{ql_invII}.
Hence, there exists $a_1b_1\in U$ for $a_1\in A-A_0$ and $b_1\in B$.
Write
$$
a_1=\left(
\begin{tabular}{cc}
$1$ & $\tilde{a}_1$\\
$0$ & $1$
\end{tabular}
\right)
$$
for some $\tilde{a}_1\in\hbox{End}(W,V)-\hbox{End}(W,\left<v_1\right>)$.
As above, 
$$
\left(
\begin{tabular}{cc}
$1$ & $\tilde{a}_1g-\tilde{a}_1$\\
$0$ & $1$
\end{tabular}
\right)\in U
$$
for all $g\in \hbox{SL}(W)$. One can choose
$g$ such that $\tilde{a}_1 g-\tilde{a}_1\notin \hbox{End}(W,\left<v_1\right>)$.
This shows that $U\cap A=A$. It follows that for any $x\notin V_\mathbb{R}$ and any $v\in V_\mathbb{R}$,
there exists $u\in U_\mathbb{R}$ such that $ux=x+v$.

Finally, we finish the proof in the Case 2.
By Lemma \ref{ql_invII}, 
$$
\mathcal{ V}=\left< \alpha \ell_2+\beta \ell_3,\ell_4\right>
$$
for some $\alpha,\beta\in \mathbb{C}$.
Let $v\ne 0$ be a vector such that 
$$
\ell_2(v)=\ell_3(v)=\ell_4(v)=0.
$$
Then $v\in V\cap V^\perp$,
and since $Q$ is nondegenerate, it follows that
$W\nsubseteq \left<v\right>^\perp$.
Note that $L|_V=\ell_4|_V=0$. In particular, $L|_W\ne 0$.
Fix $(a,b)\in\mathbb{R}^2$. Take $z\in\mathcal{P}(\mathbb{Z}^4)-V$.
Let $z=z'+z''$ for $z'\in V$ and $z''\in W-\{0\}$. Since $L|_W\ne 0$, there exists
$g\in\hbox{SL}(W_\mathbb{R})$ such that $L(gz'')=b$. 
Also we can choose $g_n\in\hbox{SL}(W_\mathbb{R})$ such that $g_n\rightarrow g$
and $Q(g_nz'',v)\ne 0$.
Take $f_n\in F^o_{0\mathbb{R}}$
such that $f_n|_{W_\mathbb{R}}=g_n$ and $u_n\in U_\mathbb{R}$ such that
$$
u_n(f_nz)=f_nz+t_nv
$$
where $t_n=\frac{a-Q(f_nz)}{2Q(g_nz'',v)}$.
Then 
\begin{eqnarray*}
Q(u_nf_nz)&=&Q(f_nz+t_nv)=Q(f_nz)+2t_nQ(f_nz,v)\\
&=&Q(f_nz)+2t_nQ(g_nz'',v)=a.
\end{eqnarray*}
Also 
$$
L(u_nf_nz)=L(f_nz'+f_nz''+t_nv)=L(g_nz'')\rightarrow b.
$$
This shows (\ref{ql_eqII}).

\underline{{\sc Case 3:} $\dim \mathcal{ V}=3$.}
Let $V=\left<v\right>$. By Lemma \ref{ql_invII}, $\mathcal{ V}=\left<\ell_2,\ell_3,\ell_4\right>$.
It follows that $Q(v)=0$. Since $F$ has no nontrivial characters, $Fv=v$. 
We show that $W$ is a $\mathbb{Q}$-subspace.
For any $\sigma\in\hbox{Gal}(\mathbb{C},\mathbb{Q})$, $W^\sigma$ is 
$F_0$-invariant (because $F_0$ is connected). Since $W\cap W^\sigma\ne \{0\}$ and
$W$ is $F_0$-irreducible, it follows that $W=W^\sigma$.
Thus, $W$ is defined over $\mathbb{Q}$.

For any $z\in W_\mathbb{Z}-\{0\}$, define
$$
F_0^z=\{g\in F^o_{0\mathbb{R}}: Q(gz,v)\ne 0\}.
$$
We claim that $F_0^z$ is not empty.
Suppose that $Q(gz,v)=0$ for every $g\in F^o_{0\mathbb{R}}$. 
Since $F^o_{0\mathbb{R}}$ is Zariski dense in $F_0$, $\left<F^o_{0\mathbb{R}}z\right>=W_\mathbb{R}$.
Thus, $W_\mathbb{R}\subseteq \left< v\right>^\perp$. Since $Q(v)=0$, the vector $v$ lies in the radical of $Q$.
This gives a contradiction. Thus, $F_0^z$ is a nonempty Zariski open
subset of $F^o_{0\mathbb{R}}$. It follows that $F_0^z$ is dense in $F^o_{0\mathbb{R}}$ in Euclidean topology.
Put 
$$
F_0^\infty=\bigcap \left\{F_0^z:z\in \mathcal{P}(W_\mathbb{Z})\right\}.
$$
This set is dense in $F^o_{0\mathbb{R}}$ by Baire Category Theorem.

Consider the orbit $O=[L]^{F^o_{0\mathbb{R}}}$ in the real projective space $\mathbb{P}(W^*_\mathbb{R})$ where
$W^*_\mathbb{R}$ is the dual space of $W_\mathbb{R}$. 
If $O$ consists of rational points, then $[L]^{F^o_{0\mathbb{R}}}=[L]$
and $[L]^{F_{0}}=[L]$, which is a contradiction.
Thus, there exists $g\in F^o_{0\mathbb{R}}$ such that the linear form $L(gx)$, $x\in W$, is not
proportional to a rational form. Then $L(g\mathcal{P}(W_\mathbb{Z}))$ is dense
in $\mathbb{R}$. Therefore, $\overline{L(F^o_{0\mathbb{R}}\mathcal{P}(W_\mathbb{Z}))}=\mathbb{R}$
and $\overline{L(F_0^\infty \mathcal{P}(W_\mathbb{Z}))}=\mathbb{R}$.

Every element $u\in U$ acts on $\mathbb{C}^4$ as follows:
$uv=v$ and for $w\in W$, $uw=w+l_u(w)v$ for some linear form $l_u$ on $W$.
Using this notations, the action of $F_0$ on $U$ by conjugation corresponds to the
usual action on the space of linear forms spanned by $l_u$, $u\in U$.
Thus, this action is irreducible.
Note that $U\ne 1$ because $\mathcal{ L}^U=\mathcal{ V}$.
It follows that for every linear form $l$ on $W$, there exists $u\in U$ such that
$uw=w+l(w)v$ for $w\in W$. In particular,
for any $w\in W$ and $t\in\mathbb{R}$, there exists $u\in U_\mathbb{R}$
such that $uw=w+tv$.

Now we are ready to finish the proof of Case 3.
Fix $(a,b)\in\mathbb{R}^2$. There exist $g_n\in F_0^\infty$ and $z_n\in \mathcal{P}(W_\mathbb{Z})$
such that $L(g_nz_n)\rightarrow b$. Take $u_n\in U_\mathbb{R}$ such that 
$u_ng_nz_n=g_nz_n+t_nv$ where $t_n=\frac{a-Q(g_nz_n)}{2Q(g_nz_n,v)}$.
Then 
$$
Q(u_ng_nz_n)=Q(g_nz_n+t_nv)=Q(g_nz_n)+2t_nQ(g_nz_n,v)=a,
$$
and 
$$L(u_ng_nz_n)=L(g_nz_n)\rightarrow b.
$$
This shows (\ref{ql_eqII}).

\underline{{\sc Case 4:} $\dim \mathcal{ V}=4$.}
By Lemma \ref{ql_ss}, $F$ is semisimple. From Table 1,
$F$ is one of the types $A_1$, $A_1+A_1$, $C_2$.
The first case is impossible because $F$ contains a $2$-dimensional
unipotent subgroup. Denote by $\mathfrak{f}$ the Lie algebra of $F$.

Let $F$ be of type $C_2$. Then 
$\mathfrak{f}=g\mathfrak{sp}(4,\mathbb{C})g^{-1}$ for some
$g\in\hbox{SL}(4,\mathbb{C})$. Equivalently,
$\mathfrak{f}=\mathfrak{sp}(S,\mathbb{C})$ for a
nondegenerate symplectic form $S$ over $\mathbb{C}$. The form $S$ is
proportional to a real form because $F$
is defined over $\mathbb{R}$.
Fix $(a,b)\in\mathbb{R}^2$. By Lemma \ref{ql_im}, there exists $x\in\mathbb{R}^4$
such that $Q(x)=a$ and $L(x)=b$. Take $x_n\in\mathbb{Q}^4-\{0\}$ such that $x_n\rightarrow x$,
and $y_n\in\mathbb{Q}^4-\left<x_n\right>$ such that $S|_{\left<x_n,y_n\right>}$ is
nondegenerate. Let $V_n=\left<x_n,y_n\right>$. Take $z_n\in\mathcal{P}(V_{n\mathbb{Z}})$.
Every element of the form $g\oplus id$ for $g\in\hbox{SL}(V_{n\mathbb{R}})$
with respect to decomposition $V_n\oplus V_n^\perp$ is in $\hbox{Sp}(S,\mathbb{R})$
for any $g\in\hbox{SL}(V_{n\mathbb{R}})$. Thus, there exists $f_n\in F^o_\mathbb{R}=\hbox{Sp}(S,\mathbb{R})$ such that
$x_n=f_nz_n$. Then $Q(f_nz_n)\rightarrow a$ and $L(f_nz_n)\rightarrow b$. This shows (\ref{ql_eqII}).

Let $F$ be of type $A_1+A_1$. 
Using an argument as in Proposition \ref{ql_subgr},
$F=\hbox{SO}(S)$ for some real nondegenerate quadratic form $S$. 

Let $\mathfrak h$ be the Lie algebra of $H$. It is easy to check that the normalizer of $\mathfrak h$
in $\mathfrak{sl}(4,\mathbb{C})$ is
\begin{equation}\label{eq_norm}
\left\{\left(
\begin{tabular}{cccc}
$u$ & $x$ & $y$ & $t$\\
$0$ & $v$ & $0$ & $y$\\
$0$ & $0$ & $-v$ & $x$\\
$0$ & $0$ & $0$ & $-u$
\end{tabular}
\right):u,v,x,y,t\in\mathbb{C}\right\},
\end{equation}
and the centralizer of $\mathfrak h$ in $\mathfrak{sl}(4,\mathbb{C})$ is
\begin{equation}\label{eq_cent}
\left\{\left(
\begin{tabular}{cccc}
$0$ & $x$ & $y$ & $t$\\
$0$ & $0$ & $0$ & $y$\\
$0$ & $0$ & $0$ & $x$\\
$0$ & $0$ & $0$ & $0$
\end{tabular}
\right):x,y,t\in\mathbb{C}\right\}.
\end{equation}

Since $F_\mathbb{R}$ contains closed unipotent subgroup $g^{-1}H_\mathbb{R}g$, it cannot be compact.
Thus, $F_\mathbb{R}$ is isomorphic to either $\hbox{SO}(3,1)_\mathbb{R}$ or $\hbox{SO}(2,2)_\mathbb{R}$. Recall that
$\hbox{SO}(3,1)_\mathbb{R}$ is isogenous to $\hbox{SL}(2,\mathbb{C})$, and
$\hbox{SO}(2,2)_\mathbb{R}$ is isogenous to $\hbox{SL}(2,\mathbb{R})\times\hbox{SL}(2,\mathbb{R})$.
In both cases, the group $g^{-1}H_\mathbb{R}g$ is a maximal unipotent subgroup of $F$.
Let $\mathfrak{h}_\mathbb{R}$ be the Lie algebra of $H_\mathbb{R}$.
Consider the map 
$$
\phi:N_{SL(4,\mathbb{R})}(gH_\mathbb{R}g^{-1})^o\rightarrow \hbox{GL}(\hbox{Ad}(g)\mathfrak{h}_\mathbb{R})
$$
defined by $\phi(n)h=\hbox{Ad}(n)h$. By (\ref{eq_norm}) and (\ref{eq_cent}), the image of $\phi$
is isogenous to $\mathbb{R}_{>0}\times\mathbb{R}_{>0}$. If $F_\mathbb{R}$ were isogenous to $\hbox{SL}(2,\mathbb{C})$, then
$\phi (N_{F_\mathbb{R}}(gH_\mathbb{R}g^{-1}))^o$ would be isogenous to $\mathbb{C}^\times\simeq \hbox{SO}(2,\mathbb{R})\times \mathbb{R}_{> 0}$.
To check the last statement, one can note that it is obvious when $g^{-1}H_\mathbb{R}g$ is the subgroup of upper unipotent matrices
in $\hbox{SL}(2,\mathbb{C})$, and  $g^{-1}H_\mathbb{R}g$ is conjugate to this subgroup.
This shows that $S$ has signature $(2,2)$.

Denote $E=\hbox{SO}(Q_0)$ and $G=\hbox{SL}(4,\mathbb{C})$.
Then $F=g^{-1}g_1^{-1}E g_1 g$ for some $g_1\in G$.
First, we show that $g_1\in EC_G(H)$. We have $H$ and $g_1Hg_1^{-1}\subseteq E$.
Since maximal unipotent subgroups are
conjugate, $eHe^{-1}=g_1Hg_1^{-1}$ for some $e\in E$. Thus, $e^{-1}g_1\in N_G(H)$,
and $g_1\in E N_G(H)$. Without loss of generality, $g_1\in N_G(H)$.
Consider the map $\psi:N_G(H)\rightarrow \hbox{GL}(\mathfrak{h})$ defined
by $\psi(n)h=\hbox{Ad}(n)h$. Since this map is
algebraic, its image is an algebraic  subgroup in $\hbox{GL}(\mathfrak{h})$.
By (\ref{eq_norm}) and (\ref{eq_cent}), $\dim \psi (N_G(H))=2$ and $\psi (N_G(H))^\circ$
is generated by diagonal matrices. Thus, $\psi (N_G(H))^\circ$ is abelian.
Let
$$
T=\{\hbox{diag}(u,v,v^{-1},u^{-1}):u,v\in\mathbb{C}^\times\}.
$$
This is a maximal torus of $E$. Since $\psi(T)$ has dimension $2$, it is a maximal
torus of $\hbox{GL}(\mathfrak h)$, and so is $\psi (g_1^{-1}Tg_1)$.
Since they commute, 
$$
\psi (g_1)^{-1}\psi (T)\psi(g_1)=\psi (T),
$$
i.e. $\psi (g_1)$ normalizes $\psi (T)$. The normalizer of $\psi (T)$
in $\hbox{GL}(\mathfrak{h})$ is generated by $\psi (T)$ and the transformation
that permutes two elements of the basis of $\mathfrak h$.
It is easy to see that this transformation is in $\psi (E)$. Thus,
$\psi (g_1)\in \psi (E)$. Since $\ker \psi=C_G(H)$, it follows that
$g_1\in EC_G(H)$.

It follows from Lemma \ref{ql_invII} that $H$ has a unique fixed vector $v$
(up to a scalar multiple). Then for $c\in C_G (H)$,  $cv=\lambda v$
for some $\lambda\in \mathbb{C}$. Suppose that $\mu\in\mathbb{C}-\{\lambda\}$
be an eigenvalue of $c$. The complex eigenspace corresponding to
$\mu$ and $\lambda$ are $H$-invariant.
Each of these subspaces contains a nonzero vector fixed by $H$ (because $H$ is unipotent).
This contradicts Lemma \ref{ql_invII}.
Thus, $c$ has a unique eigenvalue, and $c=(\lambda I)c_0$ where
$I$ is the identity matrix, and $c_0\in C_G(H)$ is unipotent.
Let $C_0$ be the set of unipotent elements of $C_G(H)$. If $c_1,c_2\in C_0$,
$c_1v=c_2v=v$, and $c_1c_2^{-1}v=v$. Then $c_1c_2^{-1}$ is unipotent.
Thus, $C_0$ is a subgroup. We have $C_G(H)=Z(G)C_0$. Clearly, $C_0$
is unipotent algebraic subgroup, so that it is connected. By (\ref{eq_cent}), 
$C_0=Hu_\mathbb{C}$ where
$$
u_t=\left(
\begin{tabular}{cccc}
$1$ & $0$ & $0$ & $t$\\
$0$ & $1$ & $0$ & $0$\\
$0$ & $0$ & $1$ & $0$\\
$0$ & $0$ & $0$ & $1$
\end{tabular}
\right).
$$
Thus, $g_1\in EZ(G)u_\mathbb{C}$. Then
\begin{eqnarray*}
F&=&g^{-1}g_1^{-1}Eg_1g=g^{-1}u_t^{-1}\hbox{SO}(Q_0)u_tg=g^{-1}\hbox{SO}(Q_0+tx_4^2)g\\
&=&\hbox{SO}(Q+tL^2)
\end{eqnarray*}
for some $t\in\mathbb{C}$.
Since $F$ is a $\mathbb{Q}$-group, the quadratic form $Q+tL^2$ is proportional
to a rational form. This is a contradiction.
The proposition is proved.
\end{proof}

\section{Conclusion}

Let $(Q,L)$ be a pair such that $Q$ is a nondegenerate quadratic form, and
$L$ is a nonzero linear form. It would be interesting to obtain neccesary and
sufficient conditions for the property
$$
\overline{\{(Q(x),L(x)):x\in\mathcal{P}(\mathbb{Z}^d)\}}=\mathbb{R}^2
$$
to hold. In this context, we formulate the following conjecture:
\begin{con} \label{ql_con}
Let $Q$ be a nondegenerate quadratic form, and $L$ a nonzero linear form in
dimension $d\ge 4$. Suppose that 
\begin{enumerate}
\item[1.] For every $\beta\in\mathbb{R}$, $Q+\beta L^2$ is indefinite.
\item[2.] For every $(\alpha,\beta)\ne (0,0)$, with $\alpha,\beta\in\mathbb{R}$, $\alpha Q+\beta L^2$ is not rational.
\end{enumerate}
Then $\{(Q(x),L(x)):x\in\mathcal{P}(\mathbb{Z}^d)\}$ is dense in $\mathbb{R}^2$.
\end{con}

The first condition in the conjecture is neccesary for the density to hold.
Indeed, suppose that $Q+\beta L^2$ is definite (say, positive definite) for some $\beta\in\mathbb{R}$.
By Proposition \ref{ql_can},
\begin{equation}\label{eq_last0}
Q=\ell_1^2+\ldots+\ell_{d-1}^2-\beta\ell_d^2\;\;\hbox{and}\;\;L=\ell_d
\end{equation}
for some linearly independent linear forms $\ell_i$, $i=1,\ldots,d$.
If $\bar{0}\in\mathbb{R}^2$ is an accumulation point of the set $\{(Q(x),L(x)):x\in \mathbb{Z}^d\}$,
then $\bar{0}\in\mathbb{R}^d$ is an accumulation point of the set 
$\{(\ell_1(x),\ldots,\ell_d(x)):x\in\mathbb{Z}^d\}$ which is impossible.

Due to Theorem \ref{ql_m}, it remains to prove Conjecture \ref{ql_con} in the case when
$(Q,L)$ is of type (II), and $Q|_{L=0}$ is positive definite, i.e.
\begin{equation}\label{eq_last}
Q=\ell_1^2+\ldots+\ell_{d-2}^2+\ell_{d-1}\ell_d\;\;\hbox{and}\;\;L=\ell_d
\end{equation}
for some linearly independent linear forms $\ell_i$, $i=1,\ldots,d$.
The method of the proof of Proposition \ref{ql_II} with minor modifications 
allows to prove Conjecture \ref{ql_con} in dimension $4$. However,
one should keep in mind that the method of reduction to lower dimension (Proposition \ref{ql_red})
fails to work in this case. Namely, for every $d\ge 4$, there exist pairs $(Q,L)$
of dimension $d$ that satisfy the conditions of
Conjecture \ref{ql_con}, but for every rational subspace $V$ of codimension $1$, the set
$\{(Q(x),L(x)):x\in V_\mathbb{Z}\}$ is not dense in $\mathbb{R}^2$.
To construct such an example, we take a pair $(Q,L)$ as in (\ref{eq_last}) such 
that the space $\left<\ell_1,\ldots,\ell_{d-2},\ell_d\right>$ contains no nonzero rational forms,
and $(Q,L)$ satisfies the second condition of Conjecture \ref{ql_con}. Clearly, such pairs exist.
In fact, such pairs are generic in suitable sense. Let $\ell$ be nonzero rational form.
It is easy to see from Lemma \ref{ql_red} that the pair $(Q|_{\ell=0},L|_{\ell=0})$ is
of type (II) iff $Q|_W$ is degenerate where $W$ is the space $\{\ell=0,L=0\}$.
The forms $\ell_i|_W$, $i=1,\ldots,d-2$, are linearly
independent because $\ell\notin\left<\ell_1,\ldots,\ell_{d-2},\ell_d\right>$.
Therefore, 
$$
Q|_W=(\ell_1^2+\ldots+\ell_{d-2}^2)|_W
$$
is nondegenerate, and $(Q|_{\ell=0},L|_{\ell=0})$ is of type (I).
Since $Q|_{L=0}$ is positive definite, the pair $(Q|_{\ell=0},L|_{\ell=0})$ is as in (\ref{eq_last0}).
As we saw above, density fails in this case.

\end{document}